\documentclass[10pt,review]{elsarticle}
\usepackage{lineno,hyperref,url}
\usepackage[hmargin=1.3cm]{geometry}
\usepackage{setspace}
\usepackage{rotating}
\usepackage{amssymb}
\usepackage{amsmath,amsthm}
\usepackage{amsfonts}
\usepackage{graphicx}
\usepackage{graphics}
\usepackage{bm}
\usepackage{color}
\usepackage [latin1]{inputenc}
\usepackage{todonotes}
\usepackage[format=plain,justification=raggedright,singlelinecheck=false]{caption}
\usepackage[justification=centering]{caption}
\makeatother \numberwithin{equation}{section}

 \modulolinenumbers[5]
 \journal{Elsevier}
\begin{document}

\begin{frontmatter}

\title{A Note on the Direct Approximation of Derivatives in \\Rational Radial Basis Functions Partition of Unity Method}

\author[mymainaddress]{Vahid Mohammadi\corref{mycorrespondingauthor}}
\cortext[mycorrespondingauthor]{Corresponding author}
\ead{vahid.mohammadi@sru.ac.ir}

\author[mymainaddress-3]{Stefano De~Marchi}
\ead{stefano.demarchi@unipd.it}

\address[mymainaddress]{Department of Mathematics, Faculty of Science,
Shahid Rajaee Teacher Training University, Tehran, 16785-163, Iran}

\address[mymainaddress-3]{Department of Mathematics "Tullio Levi-Civita", University of Padua, Italy}

\begin{abstract}
This paper proposes a Direct Rational Radial Basis Functions Partition of Unity (D-RRBF-PU) approach to compute derivatives of functions with steep gradients or discontinuities. The novelty of the method concerns how derivatives are approximated. More precisely, all derivatives of the partition of unity weight functions are eliminated while we compute the derivatives
of the local rational approximants in each patch.
As a result, approximate derivatives are obtained more easily and quickly than those obtained in the standard formulation. The corresponding error bounds are briefly discussed. Some numerical results are presented to show the technique's potential.

\end{abstract}

\begin{keyword}
 Rational RBF-PU method \sep Direct RRBF-PU approximation.
\MSC[2020]  32E30.
\end{keyword}
\end{frontmatter}

\section{Introduction}
The RBF method has been widely used to interpolate scattered data from high-dimensional spaces. The rational approximation is more effective and applicable when approximated with univariate functions with steep gradients or singularities. As detailed in \cite{Buhmannetal,DeMarchi,Jakobson}, a possible choice is to use rational interpolation, where the numerator and the denominator are polynomials. Some approaches are available to find the unknown coefficients related to the polynomials. One is
the Padé approximation, and the other is a least-squares technique \cite{Pade}. Of course, it is more difficult to implement in high-dimensional spaces because it depends on mesh generation (cf. \cite{DeMarchi}).
In \cite{Sarra}, a simplified RRBFs-PU method has been proposed, where the constant PU weight functions are considered to approximate the derivatives.
The authors of \cite{DeMarchi} have extended and analyzed the RRBFs-PU method, which is mesh-independent. It depends only on the number and location of scattered data distributed in the computational domain. Accordingly, it could be executed in high-dimensional spaces. This technique is also extended to conditionally positive definite (CPD) functions, i.e., polyharmonic spline (PHS) radial kernels augmented by polynomials \cite{Mir-Mirzaei}.
In our investigation, we rely on \cite{DeMarchi} and \cite{Mirzaei-D-RBF-PU}, but we apply another approach that allows us to describe a new formulation for the RRBF-PU method, the Direct RRBF-PU, shortly DRRBF-PU.
To find the approximate value of a derivative in the RRBF-PU, we use a direct form (cf. \cite{Mirzaei-D-RBF-PU})) instead of its standard version. This reduces the computational cost compared to that of the RRBF-PU method, and it will be observed that the error is at least of the same order as its worst local approximant.
It is worth mentioning that the authors in \cite{Liu} have recently found a compact RBF-PU technique to solve parabolic equations on surfaces. We also cite the paper \cite{DeMarchi1} that has introduced a method for approximating discontinuous functions by variably scaled discontinuous kernels (VSDK).

The remainder of the paper is presented as follows. In Section \ref{Sec-0}, an RRBF-PU approach is given. The main idea, that is, the DRRBF-PU technique is proposed in Section \ref{Sec-1}. Also, a brief discussion related to the error bounds is proposed in this section. In the next part, two examples are presented in two dimensions. Finally, a conclusion is given in Section \ref{Sec-C}.

\section{The RRBF-PU method}\label{Sec-0}
This section briefly reviews the RRBFs-PU method proposed in \cite{DeMarchi}.

We assume that the computational domain, say $\Omega$, is partitioned into the overlapping subdomains, that is, $\Omega \subseteq \bigcup_{l=1}^{N_{c}}\Omega_{l}$, with $N_c$ the number of these subdomains. Also,
assume that a family of nonnegative weight functions such as $\omega_{\ell} \in C^{k}(\mathbb{R}^{d})$ forms a partition of unity (PU), where
$supp(\omega_{l}) \subseteq \Omega_{l}$, which are open, bounded and regular (cf. \cite[Definition 15.18]{wendland2004scattered} or \cite[Chapter 29]{fasshauer2007meshfree}). Moreover,
we consider $\{\omega_{\ell}\}_{l=1}^{N_{c}}$ is {\it $k$-stable PU} w.r.t. $\{\Omega_{\ell}\}$. That is,
 for every multi-index $\alpha \in \mathbb{N}^{d}_{0}$, where $|\alpha|
\leq k$, there exists a constant $C_{\alpha}$ such that
\begin{equation}\label{weight-PU}
||D^{\alpha}\omega_{\ell}||_{L_{\infty}(\Omega)} \leq C_{\alpha}
\rho^{-|\alpha|}_{\ell},
\end{equation}
in which $\rho_{\ell}:=\frac{1}{2}{\sup _{x,y \in
        \Omega_{\ell}}}{\left\| {x - y} \right\|_2}$
\cite{fasshauer2007meshfree,Mirzaei-D-RBF-PU,wendland2004scattered}.

To construct the approximation, a set of trial points
$X=\{\bm{x}_{1},\bm{x}_{2},...,\bm{x}_{N}\} \subset \Omega \subset
\mathbb{R}^{d}$ is assumed.

Also, the set of trial points falling in
the patch $\Omega_{\ell}$ is denoted by $X_{\ell}=X \cap
\Omega_{\ell}$. Accordingly, the index set corresponding to $\Omega_{\ell}$ for $1 \leq \ell \leq N_{c}$ is defined to be $J_{\ell}:=\{j \in \{1,2,...,N\}:\bm{x}_{j} \in X_{\ell}\}$ (cf.
\cite{Mirzaei-D-RBF-PU}).

Now, we are ready to provide the RRBF-PU approximation for
a real function, $f \in
C^{k}(\Omega)$, $k$ being the order of smoothness (note that the underlying function might not be also not continuous) \cite{DeMarchi}:
\begin{equation}\label{1}
f(\bm{x}) \approx S_{f,X}(\bm{x})=\displaystyle \sum_{l=1}^{N_{c}}\omega_{l}(\bm{x})\mathcal{R}_{l}(\bm{x}),\quad \bm{x} \in \Omega,
\end{equation}
in which $R_{l}(\bm{x})$ is the local approximation in patch $l$, and its representation is
\begin{equation}\label{2}
R_{l}(\bm{x})=\dfrac{\displaystyle \sum_{j \in J_{l}}\alpha^{l}_{j}\phi(||\bm{x}-\bm{x}_{j}||_{2})}{\displaystyle \sum_{k \in J_{l}}\beta^{k}_{k}\phi(||\bm{x}-\bm{x}_{k}||_{2})}
\end{equation}
In addition, in (\ref{1}), $\omega_{l}$ is considered to be a compactly
supported, nonnegative and non-vanishing weight function on $\Omega_{l}$. As the weight function, we may consider the
Shepard function
\cite{fasshauer2007meshfree,Mirzaei-D-RBF-PU,wendland2004scattered}:
\begin{equation}\label{Local-approximation-4}
\omega_{l}(\bm{x}):=\dfrac{\psi_{l}(\bm{x})}{\displaystyle\sum_{j=1}^{N_{c}}\psi_{j}(\bm{x})},
\end{equation}
in which $\psi_{l}$ is a compactly supported RBF and we selected one of the Wendland's family (\cite[Chapter 9]{wendland2004scattered}).

Furthermore, as explained in \cite{DeMarchi}, in each patch $l$, let $\mathbf{q}_{l}$ be the corresponding  eigenvector to the smallest eigenvalue of the problem $\Lambda_{l} \mathbf{q}_{l}=\lambda_{l}\Theta_{l}\mathbf{q}_{l}$ and
\begin{equation}\label{3}
{\Lambda _l} = \frac{1}{{\left\| {{\mathbf{f}_l}} \right\|_2^2}}\left( {D_{^l}^TA_{_l}^{ - 1}{D_l}} \right) + A_l^{ - 1},\,\,\,\,\,\,\,\,{\Theta _l} = \frac{1}{{\left\| {{\mathbf{f}_l}} \right\|_2^2}}\left( {D_{^l}^T{D_l}} \right) + {I_l},
\end{equation}
where the vector $\mathbf{f}_l$ is the function values at the trial points located in patch $l$. $D_{l}$ is the diagonal matrix, in which the components of $\mathbf{f}_l$ are on its diagonal. Also, $A_{l}=[\phi(||\bm{x}_{i}-\bm{x}_{j}||_{2})]_{1 \leq i,j \leq \sharp J_{l}}$ and $I_{l}$ is an $\sharp J_{l} \times \sharp J_{l}$ identity.
The Lagrange form of (\ref{2}) can be written as follows:
\begin{equation}\label{4}
R_{l}(\bm{x})=\displaystyle \sum_{j \in J_{l}}\tilde{\alpha}^{l}_{j}\phi_{R}(||\bm{x}-\bm{x}_{j}||_{2})
=\Phi^{T}_{R}(\bm{x})\mathbf{\tilde{\alpha}^{l}}
=\Phi^{T}_{R}(\bm{x})A^{-1}_{\phi_{R}}\mathbf{f}_{X_{l}},
=\Psi^{T}(\bm{x})\mathbf{f}_{X_{l}},
=\displaystyle \sum_{j \in J_{l}}\Psi_{j}(\bm{x})f_{j},
\end{equation}
in which $\Phi^{T}_{R}(\bm{x})=[\phi_{R}(||\bm{x}-\bm{x}_{j}||_{2})]
_{1 \leq j \leq \sharp J_{l}}$. We also have used $A_{\phi_{R}}\mathbf{\tilde{\alpha}^{l}}=\mathbf{f}_{X_{l}}$, and defined $\Phi^{T}_{R}(\bm{x})A^{-1}_{\phi_{R}}:=\Psi^{T}(\bm{x})$. Note that if $\phi$ is a positive definite radial function, $\phi_{R}$ is also (cf. \cite{DeMarchi}). In our implementation here, we have applied the method of diagonal increments (MDI) to bypass the ill-conditioning of the matrix $A_l$ (cf. \cite{Sarra}).
We also refer the reader to \cite{DeMarchi, Jakobson} for more details related to the RRBFs approach.

\section{The DRRBF-PU method}\label{Sec-1}
To compute the approximate derivatives of real-valued function, say $f$ at each point $\bm{x} \in \Omega$ via Eq. (\ref{1}), the standard form takes as:
\begin{equation}\label{5}
D^{\alpha}f(\bm{x}) \approx D^{\alpha}S_{f,X}(\bm{x})=\displaystyle \sum_{l=1}^{N_{c}}D^{\alpha}(\omega_{l}(\bm{x})\mathcal{R}_{l}(\bm{x})),\quad \bm{x} \in \Omega,
\end{equation}
which shows the derivative operator $D^{\alpha}$ should act on both the weight functions $\omega_{l}$ and the local approximation $\mathcal{R}_{l}$. Hence, the computational cost would be much more expensive. To avoid computing all derivatives of the PU weight
functions, we apply the same approach proposed by Mirzaei in \cite{Mirzaei-D-RBF-PU} for the non-rational case.
Thanks to that approach, the complicated calculations
of the RRBF-PU approach could be ignored when a direct scheme is applied.
That idea also reduces the computational cost
compared to the standard RBF-PU method (cf. \cite{Mirzaei-D-RBF-PU}).
Here, we use this approach for our purposes. So,  the derivatives in (\ref{1}) are approximated as
\begin{equation}\label{6}
D^{\alpha}f(\bm{x}) \approx \widehat{D^{\alpha}S_{f,X}}(\bm{x})=\displaystyle \sum_{l=1}^{N_{c}}\omega_{l}(\bm{x})D^{\alpha}\mathcal{R}_{l}(\bm{x}),\quad \bm{x} \in \Omega,
\end{equation}
where $D^{\alpha}S_{f,X}(\bm{x})\neq\widehat{D^{\alpha}S_{f,X}}(\bm{x})$ for each $\bm{x} \in \Omega$, and $D^{\alpha}\mathcal{R}_{l}$ is the local rational approximation of
$D^{\alpha}f$ in the $l$-th patch. This shows that we have directly approximated $D^{\alpha}f$ using the PU approach. To be more precise, we first find the derivatives of the local rational approximation via the RBF in each patch $\Omega_{\ell}$. After this, we blend these approximations using the PU technique through the weight functions $\omega_{l}$ to get a global approximation $\widehat{D^{\alpha}S_{f,X}}(\bm{x})$ to $D^{\alpha}f(\bm{x})$ for all $\bm{x} \in \Omega$.
From another point of view, we can only think of the action of the differential operator based on the local approximation. So, we can ignore the smooth assumption on the PU weights \cite{Mirzaei-D-RBF-PU}. It is called a {\it direct approximation} since the method applies to RRBF in the PU setting.

To discuss the error bounds for the derivatives of the function approximated via the proposed DRRBF-PU method, we can consider the following.
Assume that in each patch $\Omega \cap \Omega_{l}$, $D^{\alpha}f$ is approximated by $D^{\alpha}R_{l}$, where $\alpha \in \mathbb{N}^{d}_{0}, \; |\alpha| \leq k/2$ is a multi-index, so that
$$||D^{\alpha}f-D^{\alpha}R_{l}||_{L_{\infty}(\Omega \cap \Omega_{l})}\leq \varepsilon_{l}.$$
Now, for any $\bm{x} \in \Omega$ and $\alpha \in \mathbb{N}^{d}_{0}$ with $|\alpha| \leq k/2$, we can write
\begin{eqnarray}\label{1-th}
         |{D^\alpha }f(\bm{x}) - \widehat {{D^\alpha }{S_{f,\bm{x}}}}(\bm{x})| \le \displaystyle\sum\limits_{l = 1}^{{N_c}} {{\omega _l}(\bm{x})\left| {{D^\alpha }f(\bm{x}) - {D^\alpha }{\mathcal{R}_l}(\bm{x})} \right|}
         \le {\displaystyle\sum\limits_{l = 1}^{{N_c}} {{\omega _l}(\bm{x})\left\| {{D^\alpha }f - {D^\alpha }{\mathcal{R}_l}} \right\|} _{{L_\infty }\left( {\Omega  \cap {\Omega _l}} \right)}}\nonumber \\
         \le \mathop {\max }\limits_{1 \le l \le {N_c}} {\left\| {{D^\alpha }f - {D^\alpha }{\mathcal{R}_l}} \right\|_{{L_\infty }\left( {\Omega  \cap {\Omega _l}} \right)}}
 \le \mathop {\max }\limits_{1 \le l \le {N_c}} \varepsilon_{l},
         \end{eqnarray}
which shows that the error bound of the derivatives of the DRRBF-PU technique is at least as good as its
worst local interpolation. This bound is derived using only the PU property of the weight functions. So, discontinuous weights can also be utilized. It is also worth mentioning that the estimation for local errors $\varepsilon_l$ is not available, and it needs a more in-depth study in future work. Of course, in the case $|\alpha|=0$, these errors follow from \cite{DeMarchi}[Proposition 3.2].
However, in the next section, we report the numerical convergence rates for the function and its first derivative.

\section{Tests}\label{Sec-2}
We consider the ${\cal C}^6$ Mat{\'e}rn kernel:
\begin{equation}
\phi(r)=e^{-cr}(15+15(cr)+6(cr)^{2}+(cr)^{3}),
\end{equation}
 whose native space is $\mathcal{N}_{\phi}(\mathbb{R}^d)=H^{\nu}(\mathbb{R}^d)$, with $\nu=\frac{d+7}{2}$ (cf. \cite[\S 4.4, p. 41]{fasshauer2007meshfree}). Here we choose $d=2$ and thus $\mathcal{N}_{\phi}(\mathbb{R}^2)=H^{4.5}(\mathbb{R}^2)$. Moreover, we take $c=35$ as a shape parameter (in an empirical way since it gives more stable and accurate solutions).
In addition, as a weight function, we choose
\cite{fasshauer2007meshfree,Mirzaei-D-RBF-PU,wendland2004scattered}:
\begin{equation}\label{PU-Weight}
\psi_{\ell}=\psi(||.-\bm{x}_{\ell}||/\rho_{\ell}),
\end{equation}
with
\begin{equation}\label{Weight-Function-New}
\psi (r) = \left\{ \begin{array}{l}
{(1 - r)^{6}}(35r^2+18r+3),\quad 0 \le r \le 1,\\
0,\quad \quad \quad \quad \quad \quad \quad \quad \quad \quad \quad
\quad \quad r
> 1,
\end{array} \right. .
\end{equation}
the Wendland ${\cal C}^4$ function which is positive definite on $\mathbb{R}^{d}$ for $d \leq 3$. Moreover,
$\bm{x}_{\ell}$ and $\rho_{\ell}$ represent the center point of the patch
$\Omega_{\ell}$ and the radius of the patches, respectively.

\begin{figure}[h!]
    \begin{center}
        \includegraphics[width=14cm,height=4cm]{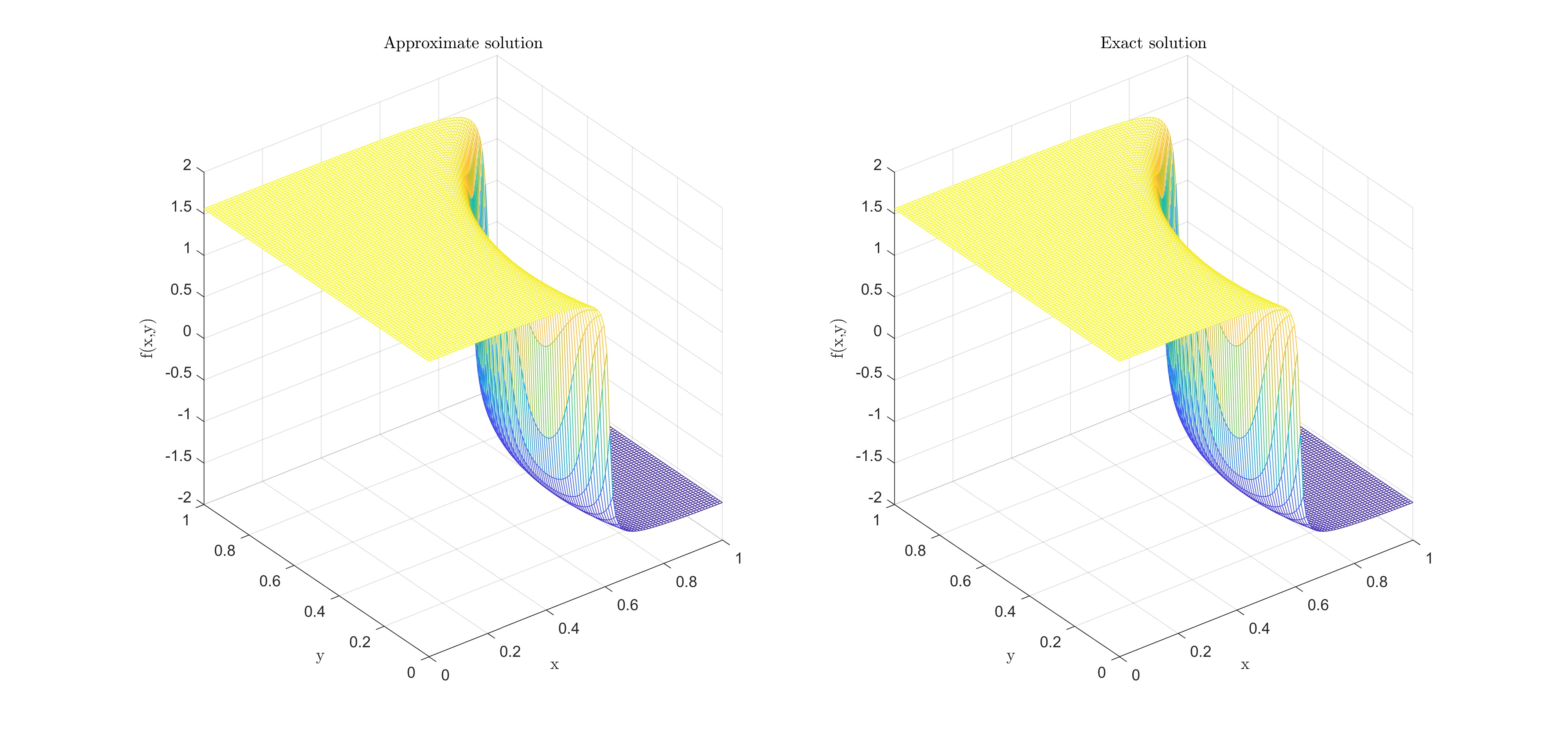}\\
        \includegraphics[width=14cm,height=4cm]{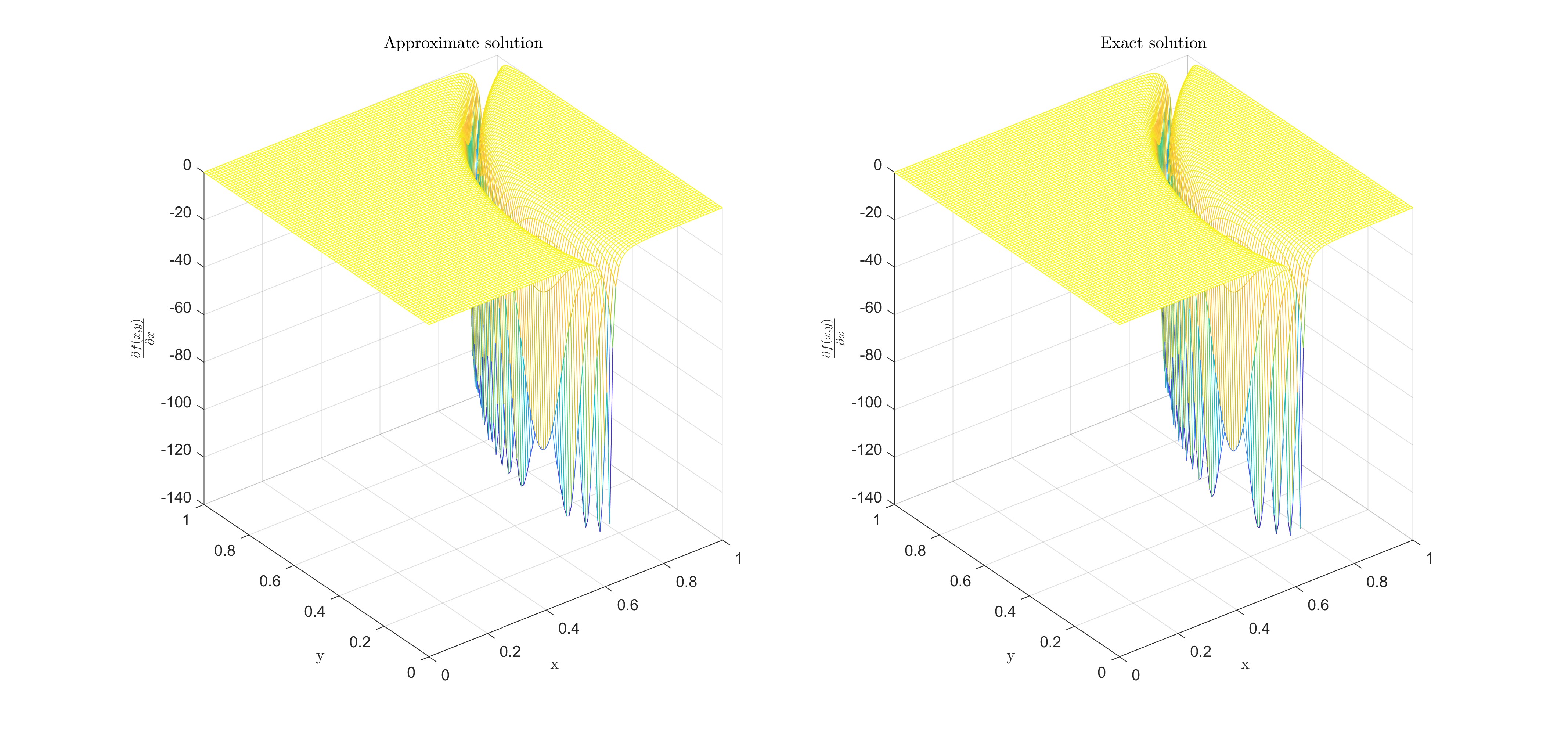}
        \caption{The approximate (left) and exact (right) solutions to function \ref{Test-1} (first row), and the approximate (left) and exact derivatives (second row). } \label{fig-1}
    \end{center}
\end{figure}

As a first test, we consider the two-dimensional function with a steep gradient \cite{Farrazandeh-Mirzaei}:
\begin{equation}\label{Test-1}
f(x,y)=\tan^{-1}\Big(125(\sqrt{(x-1.5)^2+(y-0.25)^2}-0.92)\Big),\quad (x,y) \in [0,1]^{2}.
\end{equation}
The exact and approximate solutions of $f$ and $\frac{\partial f}{\partial x}$ are plotted in Figure \ref{fig-1} using a uniform grid of $N=16641=129^2$ nodes, with to be $N_{c}=1024$ and $\rho_{\ell}=1/N_{c}$, respectively.

In addition, to measure the accuracy of the proposed technique, here we have used the following $\ell_{2}$ relative error defined by:
\begin{equation}
    ||e_{\beta}||_{2}:=
    \dfrac{||\dfrac{\partial^{\beta}f}{\partial x^{\beta}}-\dfrac{\partial^{\beta}S_{f,X}}{\partial x^{\beta}}||_{2}}{||\dfrac{\partial^{\beta}f}{\partial x^{\beta}}||_{2}},\quad \beta=0,1,
\end{equation}
so that
 $||e_{0}||_{2}$ and $||e_{1}||_{2}$ will denote the $\ell_{2}$ relative errors on the evaluation of function and its first derivative (in direction $x$), respectively. Note that the results for the first derivative in the $y$ direction are the same. So, we have not reported here for brevity. The relative errors
are computed by evaluating the numerator and the denominator on a uniform grid of $100 \times 100$ points of $[0,1]^2$.

In Table \ref{Table-1}, we report the relative errors $\ell_{2}$ of computing the approximate values of $f$, the values of the partial derivatives along $x$, and the numerical convergence rates.

\begin{table}[h!]
    \centering
    \begin{tabular}{lllllllllllll}
        \hline
        $N$&$||e_{0}||_{2}$&$Orders$&$N$&$||e_{1}||_{2}$
        & $Orders$\\
        \hline
        $1089$ &$1.56e-2$&$-$&$1089$&$2.74e-1$&$-$\\
        $4225$&$3.50e-3$&$2.10$&$4225$&$8.66e-2$&$1.66$\\
        $16641$&$1.72e-4$&$4.28$&$16641$&$2.17e-2$&$2.00$\\
        $65536$&$8.94e-6$&$4.27$&$65536$&$1.61e-3$&$3.75$\\
        \hline
    \end{tabular}
    \caption{For different $N$, $\ell_{2}$ relative errors on $f$ (left), $\ell_{2}$ relative errors on $\frac{\partial f}{\partial x}$ (right)\\ and the corresponding order of convergence.}\label{Table-1}
\end{table}

\begin{figure}[!h]
    \begin{center}
        \includegraphics[width=14cm,height=4cm]{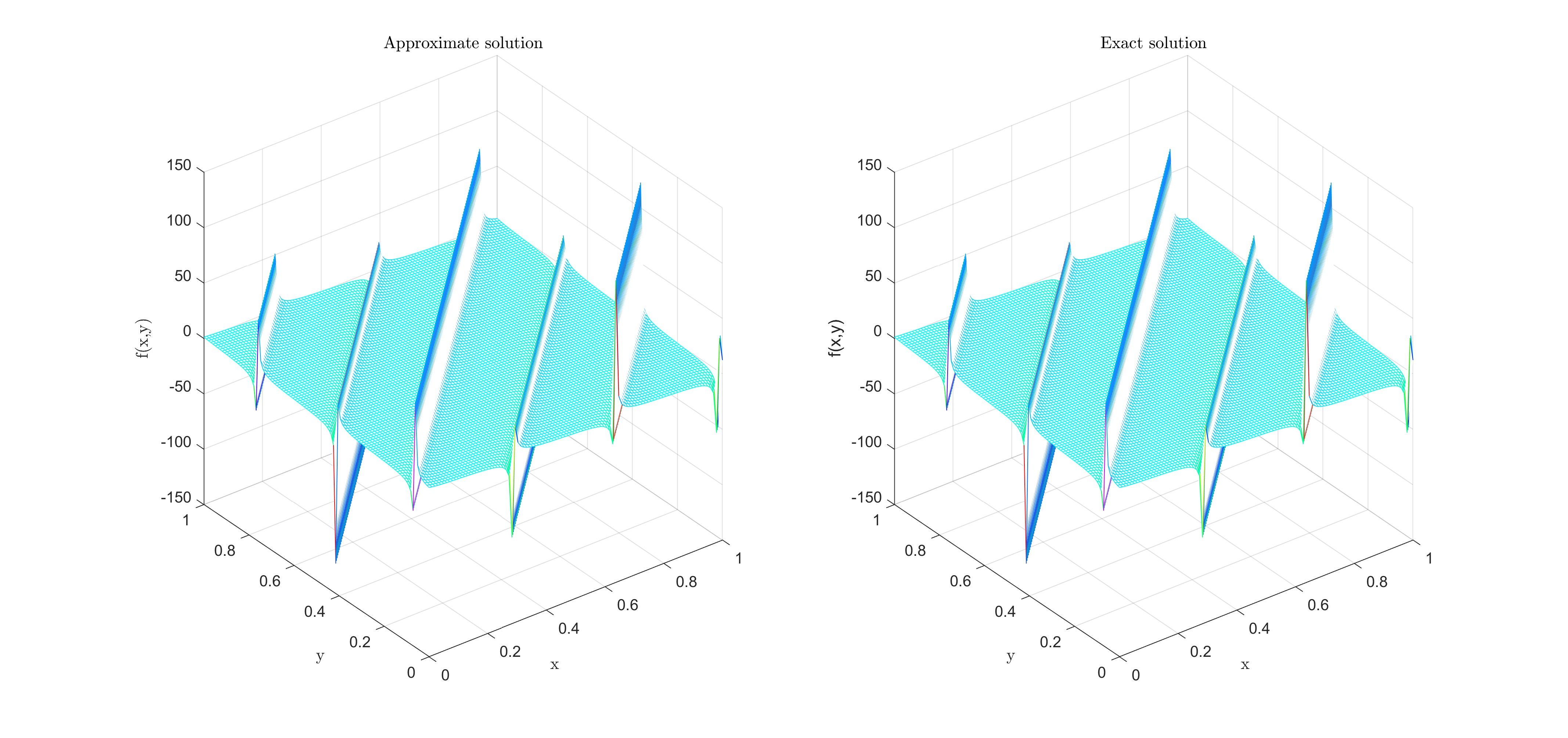}\\
        \includegraphics[width=14cm,height=4cm]{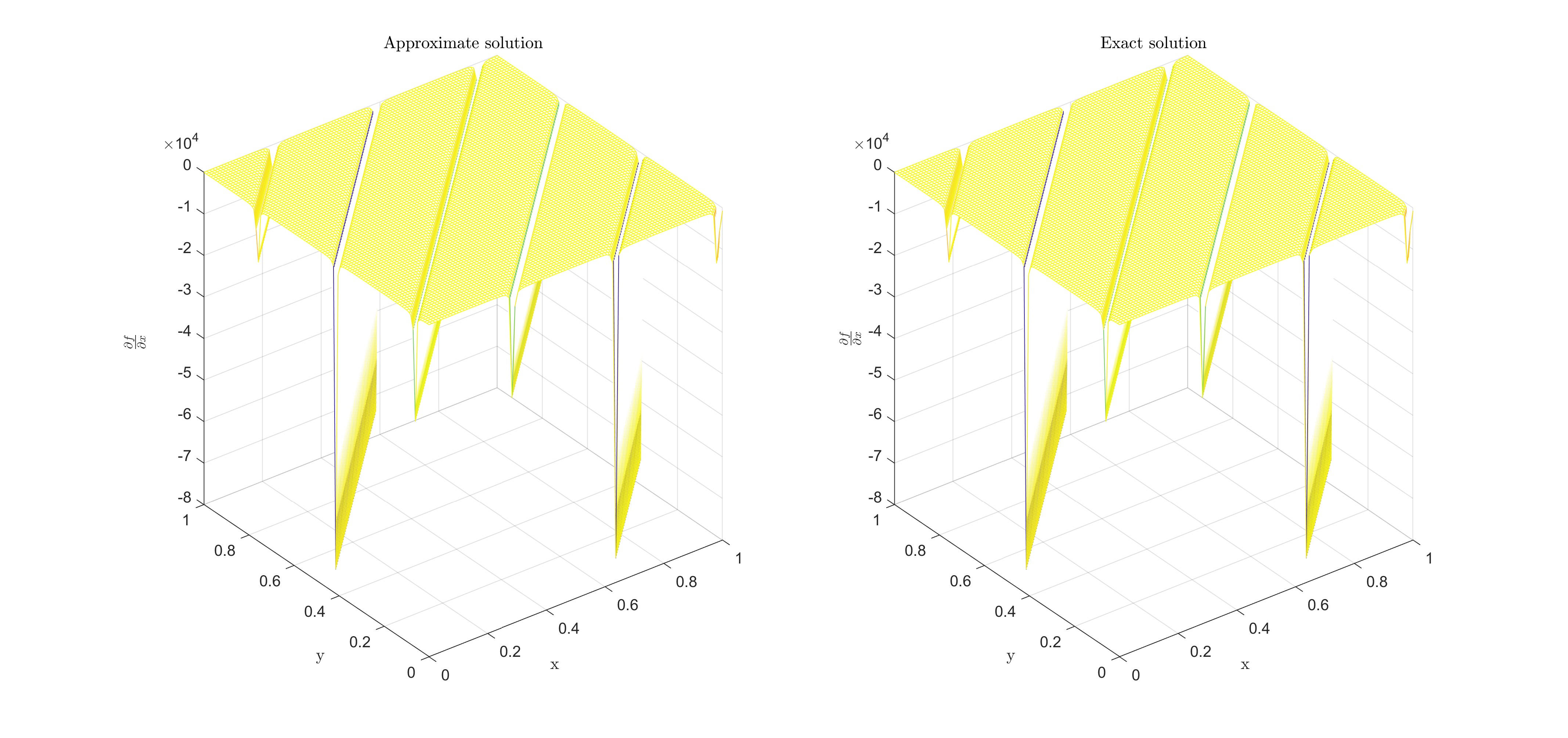}
        \caption{The approximate (left) and exact (right) solutions to function \ref{Test-2} (first row), and the approximate (left) and exact derivatives (second row). } \label{fig-2}
    \end{center}
\end{figure}

For the second test, we have considered the function also analyzed in \cite{DeMarchi}:
\begin{equation}\label{Test-2}
    f(x,y)=\frac{\tan\Big(9(y-x)+1)\Big)}{\tan (9)+1},\quad (x,y) \in [0,1]^{2}.
\end{equation}
This function has singularities across six lines $y=x+\frac{(k-9/2)\pi-1}{9}$, $k=1,\ldots,6$ (cf. \cite{Farrazandeh-Mirzaei}).
In Figure \ref{fig-2}, the exact and approximate solutions to $f$ and $\frac{\partial f}{\partial x}$ are pictured using the same $N$, $c$, $N_c$, and $\rho_{\ell}$ as in the previous example. Table \ref{Table-2} collects the corresponding $\ell_2$ errors.

All numerical tests have been produced in \verb"MATLAB" using the codes available at the link \url{https://github.com/VM-2020-MATH/DRRBFs-PU-Method}.

\begin{table}[!h]
    \centering
    \begin{tabular}{lllllllllllll}
        \hline
        $N$&$||e_{0}||_{2}$&$Orders$&$N$&$||e_{1}||_{2}$
        & $Orders$\\
        \hline
        $1089$ &$5.02e-2$&$-$&$1089$&$9.95e-2$&$-$\\
        $4225$&$4.43e-3$&$3.52$&$4225$&$9.08e-3$&$3.45$\\
        $16641$&$2.71e-4$&$4.03$&$16641$&$5.12e-4$&$4.15$\\
        $65536$&$1.43e-5$&$4.24$&$65536$&$7.24e-5$&$2.82$\\
        \hline
    \end{tabular}
     \caption{For different $N$, $\ell_{2}$ relative errors on $f$ (left), $\ell_{2}$ relative errors on $\frac{\partial f}{\partial x}$ (right)\\ and the corresponding order of convergence.}\label{Table-2}
\end{table}

\section{Conclusion}\label{Sec-C}
In this paper, we have presented the DRRBF-PU technique for computing the derivatives of functions with a steep gradient or some singularities. After outlining the numerical formulation of the method a brief discussion on the error bounds has been derived. This showed that the $\ell_2$ relative errors are of the same order as the worst local approximate constructed via the positive definite kernels. However, error estimations for approximating the derivatives regarding the fill distance are left for future work.
The two examples presented show the potential of the numerical method.
\vskip 0.2in
{\bf Acknowledgments}. The second author thanks for collaborating with the
{\it Approximation Theory and Applications} topical group of the Italian Mathematical Union, RITA the Italian Network on Approximation, and the INdAM-GNCS group.

\end{document}